\documentclass[12pt, leqno]{article}
\usepackage{amssymb,amsmath}
\usepackage[pdftex]{color}
\usepackage{latexsym}
\usepackage[pdftex]{graphicx}
\newtheorem{theorem}{Theorem}

\newtheorem{definition}{Definition}

\newtheorem{remark}{Remark}
\newtheorem{lemma}{Lemma}

\newcommand{\tx}[1]{\mbox{\;{#1}\;}} 

\newcommand{\R}{\mathbb{R}^n}
\newcommand{\rz}{\mathbb{R}}
\newcommand{\Div}{\hbox{div}\:}
\newcommand{\p}{\partial}
\numberwithin{equation}{section}
\begin{document}
\title{Isoperimetric inequalities for the principal eigenvalue of a membrane and the energy of problems with Robin boundary conditions.}
\pagestyle{myheadings}
\maketitle
\centerline{\scshape Catherine Bandle}
\medskip
{\footnotesize
\centerline{Mathematische Institut, Universit\"at Basel,}
\centerline{Rheinsprung 21, CH-4051 Basel, Switzerland}
} 
\medskip
\centerline{\scshape Alfred Wagner}
\medskip
{\footnotesize
\centerline{Institut f\"ur Mathematik, RWTH Aachen  }
\centerline{Templergraben 55, D-52062 Aachen, Germany}}
\bigskip

\abstract{An inequality for the reverse Bossel-Daners inequality is derived by means of the harmonic transplantation and the first shape derivative. This method is then applied to elliptic boundary value problems
with inhomogeneous Neumann conditions. The first variation is computed and an isoperimetric inequality is derived for the minimal energy.
\bigskip

\bigskip

\centerline {MSC2010:  49K20, 49R05, 15A42, 35J20, 35N25.}
\section{Introduction}
Bossel \cite{Bo88} and Daners \cite{Da06} extended the Rayleigh-Faber-Krahn inequality to the principal eigenvalue of the membrane with Robin boundary condition. They proved that among all domains of given volume, the first eigenvalue $\lambda$ of $\Delta \phi + \lambda \phi=0$ on $\Omega$ with $\p_\nu \phi+ \beta\phi=0$ on $\p \Omega$, where $\p_\nu$ is the outer normal derivative and $\beta$ is a positive elasticity constant, is minimal for the ball. 

Bareket \cite{Ba77} considered the eigenvalue of the same problem where $\beta$ is negative. She was able to show, that for nearly circular domains of given area the circle has the largest first eigenvalue. Recently this result was extended to higher dimensions for nearly spherical domains by Ferone, Nitsch and Trombetti \cite{ FeNiTr14}. The question whether or not the ball is optimal for all domains of the same volume remains open.

In this note we derive an isoperimetric inequality for arbitrary domains in $\R$. The proof uses the method of harmonic transplantation which is a substitute of the conformal transplantation in higher dimensions.

The method applies to a large class of variational problems related to elliptic equations with homogeneous and inhomogeneous boundary conditions. We illustrate it by means of some problems with inhomogeneous Neumann conditions. Such problems have been considered by Auchmuty \cite{Au14} in the context of  trace inequalities.

The quantities we  want to estimate in this paper are expressed by variational principles, defined for functions in $W^{1,2}(\Omega)$ and $L^2(\p \Omega)$.  For the existence of those quantities we have to require that the embedding $W^{1,2}(\Omega)$ into $L^2(\Omega)$ as well as the trace operator  $\Gamma : W^{1,2}(\Omega) \to L^2(\p \Omega)$ is compact. This is the case, cf. \cite{Au12}, when $\Omega \subset \R$ is a bounded domain such that its boundary consists of a finite number of closed Lipschitz surfaces of finite surface area. {\sl Throughout this paper we shall assume that this condition is satisfied.}

\section{Eigenvalue problem}
Consider the eigenvalue defined by
\begin{align}\label{eigenvalue}
\lambda(\Omega)= \inf_{W^{1,2}(\Omega)} \left\{\int_\Omega(|\nabla u|^2\:dx -\alpha \oint_{\p \Omega} u^2\:dS\right\} \tx{with} \int_\Omega u^2\:dx =1 \tx{and $\alpha >0$}.
\end{align}
Under our assumptions on $\Omega$ the minimum is achieved and the minimizer $u$ is a solution of
\begin{align}\label{eigenfunction}
\Delta u + \lambda(\Omega) u =0 \tx{in} \Omega, \: \p_\nu u= \alpha u \tx{on} \p \Omega.
\end{align}
If we choose a constant as a trial function in \eqref{eigenvalue}, we see immediately that $\lambda(\Omega)$ is negative.
\medskip

This problem appears in acoustics and has been discussed by M. Bareket \cite{Ba77}. She shows that for nearly circular domains obtained by surface preserving perturbations, $\lambda(\Omega)$ is largest for the circle. This result has been extended to higher dimensions in \cite{FeNiTr14}.  The main tool was the first domain variation for $\lambda(\Omega)$. 
\subsection{Domain variation and first variation}
In this section we follow closely the paper \cite{BaWa14}. Let $\Omega_t$ be a family of perturbations of the domain $\Omega\subset \R$ of the form
\begin{align}\label{yvar}
\Omega_t:=\{y=x+tv(x)+o(t):x \in \Omega,\: t\tx{small}\},
\end{align}
where $v=(v_1(x),v_2(x),\dots,v_n(x))$ is a smooth  vector field and where $o(t)$ collects all terms such that $\frac{o(t)}{t}\to 0$ as $t\to 0$. 
Note that with this notation we have
\begin{eqnarray}\label{volin}
\vert \Omega_t\vert=\vert\Omega\vert +t\:\int\limits_{\Omega}\Div\:v\:dx+o(t),
\end{eqnarray}
where $\vert\Omega\vert $ denotes the $n$ - dimensional Lebesgue measure of $\Omega$.
\newline
\newline
We say that $y: \Omega_t\to \Omega$ is volume preserving of the first order if 
\begin{align}\label{volume1}
\int_\Omega \Div v\:dx =\oint_{\p \Omega} (v\cdot \nu)=0.
\end{align}
Let $\lambda(\Omega_t)$ be the eigenvalue of a perturbed domain $\Omega_t$ (as described in \eqref{yvar} ). Let $u(t)$ be the corresponding eigenfunction. Thus $u(t)$ and  $\lambda(\Omega_t)$ solve
$$
\Delta u(t) + \lambda(\Omega_t) u(t) =0 \tx{in} \Omega_t, \: \p_{\nu_t} u(t)= \alpha\: u(t) \tx{on} \p \Omega_t,
$$
where $\nu_t$ is the outer normal of $\Omega_t$. We will use the notation $\lambda=\lambda(0)=\lambda(\Omega)$.
\newline
\newline
The first variation of $\frac{d}{dt}\lambda(\Omega_t)\big |_{t=0}=:\dot\lambda(0)$ has been computed by different authors and assumes the form
\begin{align}\label{eigenvaluev1}
\dot \lambda (0)= \oint_{\p \Omega}(|\nabla u|^2-\lambda(0)u^2-2\alpha^2 u^2-\alpha (n-1)Hu^2)(v\cdot \nu)\:dS,
\end{align}
where $H$ is the mean curvature of $\p \Omega$. From this formula we get immediately the 
\begin{lemma} \label{monotony}Let $\Omega=B_R$ be  the ball of radius $R$ centered at the origin. Suppose that
$|\Omega_t|>|B_R|$ fo small $|t|$. Then
$$
\dot{\lambda}(0)>0.
$$
In particular $\lambda(B_{R_1})> \lambda(B_{R_0})$ if $R_1>R_0$.
\end{lemma}
\noindent {\bf Proof} 
It follows from the variational characterization that the first eigenfunction is of constant sign and radial. The eigenvalue problem \eqref{eigenfunction} then reads as
\begin{eqnarray}\label{radial}
u''(r)+\frac{n-1}{r} u'(r)+ \lambda(B_R)\: u(r)=0, \quad u'(R)=\alpha u(R).
\end{eqnarray}
Moreover for the integrand in \eqref{eigenvaluev1} we have
\begin{eqnarray*}
&&\vert\nabla u\vert^2-\lambda\:u^2-2\alpha^2 u^2-\alpha (n-1)Hu^2\\
&&\qquad\qquad=
\vert u'(R)\vert^2-\lambda\:u^2(R)-2\alpha^2 \:u^2(R)-\alpha\: \frac{n-1}{R}u^2(R)\\
&&\qquad\qquad=
-\left(\lambda+\alpha^2+\alpha\: \frac{n-1}{R}\right)u^2(R)
\end{eqnarray*}
since $u'(R)=\alpha u(R)$. Thus
\begin{eqnarray*}
\dot \lambda (0)=-\left(\lambda+\alpha^2+\alpha\: \frac{n-1}{R}\right)u^2(R)\: \oint_{\p B_R}(v\cdot \nu)\:dS.
\end{eqnarray*}
Since $|\Omega_t|>|B_R|$ for small $t$, formula \eqref{volin} and then \eqref{volume1} imply
\begin{eqnarray*}
\oint_{\p B_R}(v\cdot \nu)\:dS>0.
\end{eqnarray*}
Thus 
\begin{eqnarray*}
\dot \lambda (0)>0\qquad\hbox{iff}\qquad \left(\lambda+\alpha^2+\alpha\: \frac{n-1}{R}\right)<0.
\end{eqnarray*}
This will be proved with the help of \eqref{radial}. We set $z=\frac{u_r}{u}$ and observe that
$$
\frac{dz}{dr} +z^2 +\frac{n-1}{r}z + \lambda=0  \tx{in} (0,R).
$$
At the endpoint
$$
\frac{dz}{dr}(R) +\alpha^2 +\frac{n-1}{R} \alpha  + \lambda =0.
$$ 
We know that $z(0)=0$ and $z(R)= \alpha>0$. Note that 
\begin{eqnarray}\label{zinc}
z_r(0)= -\lambda >0,
\end{eqnarray}
thus $z(r)$ increases near $0$.  Let us now determine the sign of $z_r(R)$. 
If $z_r(R)\leq 0$ then because of \eqref{zinc} there exists a number $\rho \in (0,R)$ such that $z_r(\rho) =0$, $z(\rho)>0$ and $z_{rr}(\rho) \leq 0$. 
From the equation we get $z_{rr}(\rho) =\frac{n-1}{\rho^2} z(\rho)>0$ which leads to a contradiction. Consequently
\begin{align}\label{signk}
z_r(R)= -( \alpha^2 +\frac{\alpha (n-1)}{R}+\lambda)>0.
\end{align}
Hence
$$
\dot \lambda (0)>0 
$$
for all volume increasing perturbations $\oint_{\p B_R}v\cdot \nu \:dS >0$. This completes the proof of the lemma. \hfill $\square$
\medskip

This monotonicity is opposite to the usual case where $\alpha$ is negative and it will be crucial for the upper bounds derived in the next section. 
\subsection{Harmonic transplantation and isoperimetric inequality}
In this section we recall the method of harmonic transplantation which has been deviced by Hersch\cite{He69}, (cf. also \cite{BaFl96}) to construct trial functions for variational problems of the type \eqref{eigenvalue}. To this end we need the Green's function with Dirichlet boundary condition
\begin{eqnarray}\label{Green}
G_\Omega(x,y)= \gamma( S(|x-y|)-H(x,y)),
\end{eqnarray}
where
\begin{align}
\gamma = 
\begin{cases} 
\frac{1}{2\pi}&  \tx{if} n=2\\
\frac{1}{(n-2)|\p B_1|}& \tx{if} n>2
\end{cases}
&\qquad\tx{and} \qquad S(t)= \begin{cases}
-\ln(t) &\tx{if} n=2\\
t^{2-n} & \tx{if} n>2.
\end{cases}
\end{align}
For fixed $y\in \Omega$ the funcion $H(\cdot,y)$ is harmonic.
\begin{definition} The harmonic radius  at a point $y\in \Omega$ is given by
$$
r(y)= \begin{cases}
e^{-H(y,y)} &\tx{if} n=2,\\
H(y,y)^{-\frac{1}{n-2}}&\tx{if} n>2.
\end{cases}
$$
\end{definition}
The harmonic radius vanishes on the boundary $\p\Omega$ and takes its maximum $r_\Omega$ at the harmonic center $y_h$.  It satisfies the isoperimetric inequality \cite{He69},\cite{BaFl96}
\begin{align}\label{harmonic radius}
|B_{r_\Omega}|\leq |\Omega|.
\end{align}

Note that $G_{B_R}(x,0)$ is a monotone function in $r=\vert x\vert$. Consider any radial function $\phi:B_{r_\Omega}\to \mathbb{R}$ thus $\phi(x)=\phi(r)$. Then there exists a function $\omega:\rz\to\rz$ such that
\begin{eqnarray*}
\phi(x)=\omega(G_{B_{r_\Omega}}(x,0)).
\end{eqnarray*}
To $\phi(x)$ we associate the transplanted function $U:\Omega \to \mathbb{R}$  defined  by $U(x)=\omega(G_\Omega(x,y_h))$. Then for any positive function $f(s)$, cf.. \cite{He69} or \cite{BaFl96}, the following inequalities hold true
\begin{align}\label{harmonic1}
\int_\Omega |\nabla U|^2\:dx &= \int_{B_{r_\Omega}} \vert\nabla \phi|^2\:dx\\
\int_\Omega f(U)\:dx &\geq \int_{B_{r_\Omega}} f(\phi)\:dx.\label{harmonic2}
\end{align}
For our purpose we need an estimate of $\int_\Omega f(U)\:dx$ from above. For this some auxiliary lemmata are needed.
The following notation will be used.
\begin{align*}
\Omega^t:=\{ x\in \Omega: G_\Omega(x,y_h)>t\},&\quad B^t:=\{ x\in \Omega: G_B(x,0)>t\},\\
m_\Omega(t)= |\Omega^t|, &\quad m_B(t):=|B^t|.
\end{align*}
Recall that the Green's function $G_\Omega(x,y_h)$ is  harmonic in the domain $\Omega\setminus\Omega^t$ and constant on the boundary and that $G_{B_{r_\Omega}}(x,0)$ has analogous properties. Furthermore the capacity of the two sets is given by
\begin{eqnarray*}
\rm{cap}(\Omega\setminus\Omega^t)=\frac{1}{t^2}\int\limits_{\Omega\setminus\Omega^t}\vert\nabla G_\Omega(x,y_h)\vert^2\:dx\quad\hbox{and}\quad
\rm{cap}(B_{r_\Omega}\setminus B^t)=\frac{1}{t^2}\int\limits_{B_{r_\Omega}\setminus B^t}\vert\nabla G_{B_{r_\Omega}}(x,0)\vert^2\:dx.
\end{eqnarray*}
If we use the fact that
\begin{eqnarray*}
\oint_{\partial\left(\Omega\setminus\Omega^t\right)}\partial_{\nu} G_\Omega(x,y_h)\:dS
=
\oint_{\partial\left(B_{r_\Omega}\setminus B^t\right)}\partial_{\nu} G_{B_{r_\Omega}}(x,0)\:dS
=
t,
\end{eqnarray*}
a simple computation shows that the capacities of $\Omega\setminus \Omega^t$ and $B_{r_\Omega}\setminus B^t$ are equal. Let $r_{t}$ be the radius of $B^t$, then
\begin{eqnarray}\label{cap1}
\rm{cap}(\Omega\setminus\Omega^t)=\rm{cap}(B_{r_\Omega}\setminus B^t)=
\begin{cases}
|\p B_1|\:\frac{n-2}{r_{t}^{2-n}-r_\Omega^{2-n}} &\tx{if} n>2,\\
|\p B_1|\:\left(\ln(r_{t})-\ln(r_{\Omega})\right)&\tx{if} n=2.
\end{cases}
\end{eqnarray}
The following lemma is crucial for our optimization result.
\begin{lemma} \label{cap2}
Let 
$$
\gamma := \left(\frac{|\Omega|}{|B_{r_\Omega}|}\right)^{\frac{1}{n}}.
$$
Then $m_\Omega(t) \leq \gamma ^n \:m_B(t)$ for all $t\in (0,\infty)$.
\end{lemma}

{\bf Proof} By a rearrangement argument
$$
\rm{cap}(\Omega\setminus \Omega^t) \geq \rm{cap}(B_R\setminus B_\rho),
$$
where $B_R$ is the ball with the same volume as $|\Omega |$ and $B_{\rho_{t}}$ is the ball with the same volume as $|\Omega^t|$.
From \eqref{cap1} we deduce that
\begin{align*}
\frac{1}{r_{t}^{2-n}- r_\Omega^{2-n}}\geq \frac{1}{\rho_{t}^{2-n}- R^{2-n}} &\quad\tx{if} n>2,\\
\ln(r_{t})-\ln(r_{\Omega})\geq \ln(\rho_{t})-\ln(R) &\quad\tx{if} n=2.
\end{align*}
Hence
\begin{align}\label{eq:1}
\rho_{t}^{2-n}-R^{2-n}\geq r^{2-n}_{t}-r_\Omega^{2-n}.
\end{align}
By the definitions of $R$ and $r_\Omega$
\begin{eqnarray*}
\gamma =\frac{R}{r_\Omega}\:>\:1,
\end{eqnarray*}
hence $R=\gamma \:r_\Omega$. Introducing this expression into \eqref{eq:1} we find
$$
\rho_{t}^{2-n}\geq r^{2-n}_{t}-r_\Omega^{2-n}(1-\frac{1}{\gamma^{n-2}})
\geq (\gamma\: r_{t})^{2-n}\qquad\hbox{for}\:\: n\geq 3,
$$
and
$$
\ln(r_t)\geq\ln(\rho_t)-\ln(\gamma)\qquad\hbox{for}\:\: n=2.
$$
Consequently $\rho_{t}\leq \gamma\: r_{t}$ which completes the proof. 
\hfill $\square$
\medskip

This lemma enables us to construct an upper bound for $\int_\Omega f(U)\:dx$.
\begin{lemma}
Suppose that $f(s)$ is positive and monotone increasing. Let $\phi(x):B_{r_\Omega}\to\rz$ be radial and monotone increasing. Then
\begin{align}\label{cap3}
\int_\Omega f(U)\:dx \leq \gamma^n\int_{B_{r_\Omega}} f(\phi)\:dx.
\end{align}
\end{lemma}
{\bf Proof} 
Integration along level surfaces implies
$$
\int_\Omega f(U)\:dx = -\int_0^\infty f(\omega(t))\:dm_\Omega(t)= -f(\omega(t))\:m_\Omega(t)\big |_0^\infty +\int_0^\infty f'(\omega(t))\:\omega'(t)\:m_\Omega (t)
\:dt.
$$
If $f(\omega(0))$ is bounded 
$$
\int_\Omega f(U)\:dx = f(\omega(0))|\Omega| + \int_0^\infty f'(\omega(t))\:\omega'(t)\:m_\Omega (t)\:dt.
$$
The assertion now follows from Lemma \ref{cap2}. \hfill $\square$
\medskip

\noindent We are now in position to prove
\begin{theorem}\label{estimate1}
If $\Omega\subset \R$ is any domain with maximal harmonic radius $r_\Omega$ then
$$
|\Omega|\lambda(\Omega) \leq |B_{r_\Omega}|\lambda(B_{r_\Omega}).
$$
Equality holds if and only if $\Omega$ is the ball $B_{r_\Omega}$.
\end{theorem}
\medskip 
{\bf Proof} Let $u(|x|)$ be the eigenfunction corresponding to $\lambda(B_{r_\Omega})$ and $U$ be its transplantation into $\Omega$. Then by \eqref{eigenvalue}
$$
\lambda(\Omega)\leq \frac{\int_\Omega |\nabla U|^2\:dx -\alpha \oint_{\p \Omega} U^2\:dS}{\int_\Omega U^2\:dx}.
$$ 
In view of the equality \eqref{harmonic1} the numerator becomes
$$ 
\int_{B_{r_\Omega}} |\nabla u|^2\:dx -\alpha u^2(r_\Omega)\oint_{\p \Omega}\:dS. 
$$
The isoperimetric inequality together with  \eqref{harmonic radius} implies
$$
|\p \Omega| \geq c_n |\Omega |^{\frac{n-1}{n}} \geq c_n |B_{r_\Omega}|^{\frac{n-1}{n}} = |\p B_{r_\Omega}| .
$$
From these estimates and the fact that $\lambda(B_{r_\Omega})<0$ it follows that
$$
\int_{B_{r_\Omega}} |\nabla u|^2\:dx -\alpha u^2(r_\Omega)\oint_{\p \Omega}\:dS<0.
$$
Since $u(|x|)$ is a positive radial increasing function \eqref{cap3} applies and thus
$$
\lambda(\Omega)\leq \frac{\int_{B_{r_\Omega}} |\nabla u|^2\:dx -\alpha \oint_{\p B_{r_\Omega}}u^2\:dS}{\gamma^n\int_{B_{r_\Omega}} u^2\:dx}=\gamma^{-n}\lambda(B_{r_\Omega})
$$
which completes the proof. \hfill $\square$
\newline
\newline
For the ball,  $\lambda(B_r)$ can be determined implicitly by
\begin{align}\label{ballco}
\sqrt{|\lambda(B_r)|} = \left(\alpha + \frac{n-2}{2r}\right) \frac{ I_\nu(\sqrt{|\lambda(B_r)|}r)}{I'_\nu(\sqrt{|\lambda(B_r|}r)},
\end{align}
where  $I_\nu$ denotes the modified Bessel function of order $\nu=\frac{n-2}{2}$. From Lemma  \ref{monotony} it follows that $\sqrt{|\lambda(B_r)|}$ increases as $r$ increases. Therefore
$\lim_{r\to \infty} \sqrt{|\lambda(B_r)|}r \to \infty$ as $r\to \infty$. This together with the asymptotic behavior of $I_\nu(z)$, namely $I_\nu(z) \sim \frac{e^z}{\sqrt{\pi z}}$
as $z\to \infty$, gives
$$ 
\lim_{r\to \infty} \sqrt{\vert\lambda(B_r)| }= \alpha.
$$
\begin{remark} Theorem \ref{estimate1} is weaker than the estimate 
\begin{align}\label{Bareket}
\lambda(\Omega)\leq \lambda(B_R).
\end{align}
This is a consequence of \eqref{ballco}. In fact if we set $r^n|\lambda(B_r)|=:y^2$ then
$$
y=r^{n/2}\left(\alpha +\frac{n-2}{2r}\right) \frac{I_\nu(yr^{-\nu})}{I'_\nu(yr^{-\nu})}.
$$
Since $I_\nu/I_\nu'$ is decreasing straightforward differentiation shows that $y'$ is increasing.
\end{remark}
\begin{remark} For given $|\Omega|$ it is always possible to find a domain with a large boundary surface such that $\lambda(\Omega)<\lambda(B_R)$.  This can be seen as follows. If we introduce in \eqref{eigenvalue} a constant then
$$
\lambda(\Omega) <-\frac{\alpha |\p\Omega|}{|\Omega|}.
$$
The expression at the right-hand side can be made arbitrarily small whereas $\lambda(B_R)$  is fixed for given $|\Omega|$. 
\end{remark}
\begin{remark}\label{2var} 
In \cite{BaWa14} the second variation $\ddot{\lambda}(0)$ was computed for $\alpha < 0$. In particular (7.14) applies to our problem, if we replace $\alpha$ by $-\alpha$ there.
Next we follow the arguments which led to (7.19) and obtain $\ddot{\lambda}(0)<0$.
\end{remark}
\section{Steklov type problems}
In this section we study problems with a variable weight on the boundary. Let $\rho(x)$ be
a continuous function defined in $D\supset \Omega_t$ for $|t|\leq \epsilon$. Consider boundary value problems of the type
\begin{eqnarray}\label{toreq}
\Delta u +G'(u)=0\quad\hbox{in}\:\Omega\qquad\qquad\partial_{\nu}u=\mu\:\rho(x)\quad\hbox{in}\:\partial\Omega.
\end{eqnarray}
This equation is understood in the weak sense
\begin{eqnarray}\label{torweq}
\int\limits_{\Omega}\nabla w\cdot\nabla\varphi\:dx-\int\limits_{\Omega}G'(u)\varphi\:dx
-\mu\int\limits_{\partial\Omega}\varphi\:\rho(x)\:dS=0
\end{eqnarray}
for all $\varphi\in H^{1,2}(\Omega)$. It is the Euler -Lagrange equation corresponding to the energy
\begin{eqnarray}\label{toren}
{\cal{E}}(v,\Omega):=\int\limits_{\Omega}\vert\nabla v\vert^2\:dx\:dx-2\int\limits_{\Omega}G(v)\:dx-2\mu\int\limits_{\partial\Omega}v\:\rho(x)\:dS
\end{eqnarray}
for $u\in H^{1,2}(\Omega)$. A special case where $G'(u)$ is constant appears in \cite{Au14}.  
\newline
\newline
We consider problem \eqref{toreq} in the perturbed domains $\Omega_t$ described
in \eqref{yvar}. We assume that there exists a unique solution. The corresponding energy \eqref{toren} will be denoted by ${\cal{E}}(t)$. Following \cite{BaWa14} we compute the first variation $\dot{\cal{E}}(0)$ formally.
For the calculation we refer to \cite{BaWa14} . There it has been carried out in detail
for a more general case. In particular we use Section 4.1 and formula (2.18). 
\newline
\newline
Let us decompose $v$ in its normal and tangential components $v=v^\tau +(v\cdot \nu)\nu$. Let $\Div_{\p\Omega} v= \p_iv_i -\nu_j\p_jv_i\nu_i$ be the tangential divergence. Here we use the Einstein convention. Then
\begin{eqnarray*}
v\cdot \nabla u= v^\tau\cdot \nabla^\tau u+\mu (v\cdot \nu)\rho(x).
\end{eqnarray*}
It then follows that
\begin{eqnarray}
\nonumber&&\dot{{\cal{E}}}(0)=\oint_{\p \Omega}\{|\nabla u|^2-2G(u)\}(v\cdot \nu) \:dS
-
2\oint_{\p \Omega}(v\cdot\nabla u)\:\partial_{\nu}u\:dS\\
&&\label{dotmu}\qquad-
2\mu\oint_{\p \Omega}u\:v\cdot\nabla\rho(x)\:dS
-
2\mu\int\limits_{\partial\Omega} u\:\rho(x)\Div_{\p \Omega}v^\tau\:dS \\
\nonumber&&\qquad+
2(n-1)\mu\int\limits_{\partial\Omega} (v\cdot \nu) H\:dS.
\end{eqnarray}
By  \eqref{toreq} 
\begin{eqnarray*}
-2\oint_{\p \Omega}(v\cdot\nabla u)\:\partial_{\nu}u\:dS
&=&
-2\oint_{\p \Omega}(v^{\tau}\cdot\nabla^{\tau} u)\:\partial_{\nu}u\:dS
-
2\mu\oint_{\p \Omega}(v\cdot\nu)\rho(x)\:\partial_{\nu}u\:dS\\
&=&
-2\mu\oint_{\p \Omega}(v^{\tau}\cdot\nabla^{\tau} u)\:\rho(x)\:dS
-
2\mu^2\oint_{\p \Omega}(v\cdot\nu)\rho(x)^2\:dS\\
&=&
2\mu\oint_{\p \Omega}\Div_{\partial\Omega}v^{\tau}\:u\:\rho(x)\:dS
+
2\mu\oint_{\p \Omega}(v^{\tau}\cdot\nabla^{\tau} \rho(x))\:u\:dS\\
&&-
2\mu^2\oint_{\p \Omega}(v\cdot\nu)\rho(x)^2\:dS.
\end{eqnarray*}
We compare this with \eqref{dotmu} and finally get
\begin{eqnarray}\label{var1}
\quad\dot{{\cal{E}}}(0)=\oint_{\p \Omega}\{|\nabla u|^2-2G(u)-2\mu^2\rho(x)^2-2\mu\:u\:\partial_{\nu}\rho(x)+2(n-1)\mu\:H\}(v\cdot \nu) \:dS.
\end{eqnarray}

{\sc Definition} A domain $\Omega$ is said to be {\sl critical in the class of $\Omega_t$ }  if  $\dot{{\cal{E}}}(0)=0$.
\medskip

Observe that \eqref{var1} gives a necessary condition for the solution $u$ of \eqref{toreq} in a critical domain and in particular for an extremal domain. For specific perturbation such as volume preserving perturbation the discussion of the overdetermined boundary problem is still open.
\medskip

\noindent {\sc Example} 
\smallskip

1. $B_R$ is a critical domain if the solution $u$ of \eqref{toreq} and $\rho$ are radial and if the perburtation is volume preserving in the sense of \eqref{volume1}.
\medskip

2. Let $\Omega=B_R$ and $G(u)=k u$. Then \eqref{toreq} becomes
$$
\Delta u +k =0 \tx{in} B_R, \quad \p_\nu u =\mu \rho \tx{on} \p B_R.
$$
Note that this problem has a solution if and only if $k=\mu\oint_{\p B_R} \rho\:dS=:\mu M.$  The solution is not unique and the energy is not a minimizer.
\medskip

In the next section we will investigate a relaxed formulation of a related optimization problem.

\medskip

\subsection{Isoperimetric inequalities}
In this section we reconsider the energy given in \eqref{toren}. In particular we
assume that 
$$
{\cal{E}}(\Omega) = \inf_{W^{1,2}(\Omega)}\int\limits_{\Omega}\vert\nabla v\vert^2\:dx\:dx-2\int\limits_{\Omega}G(v)\:dx-2\mu\int\limits_{\partial\Omega}v\:\rho(x)\:dS
$$
attains its minimum and that there is a unique minimizer $u$ which solves \eqref{toreq}. The aim is to derive an upper bound by means of harmonic transplantation. We shall distinguish between two cases.
\medskip

(i) $G(s)>0$.

\noindent Consider the comparison problem
\begin{align}\label{ball}
\Delta \phi +G'(\phi) =0 \tx{in} B_{r_\Omega}, \qquad\p_\nu \phi= \mu M \tx{on} \p B_{r_\Omega} \tx{where} M:=\oint_{\p B_{r_\Omega}}\rho 
\:dS.
\end{align}
Because of the extremal property of the corresponding energy,  $\phi$ is radially symmetric. The arguments of Section 2.2, in particular the inequalities \eqref{harmonic1} and \eqref{harmonic2} apply. Let $\phi(x)=\omega(G_{B_{r_\Omega}}(x,0))$ and set $U(x)=\omega(G_\Omega(x,y_h))$. Then
$$
{\cal{E}}(\Omega) \leq \int\limits_{\Omega}\vert\nabla U\vert^2\:dx-2\int\limits_{\Omega}G(U)\:dx-2\mu\int\limits_{\partial\Omega}U\:\rho(x)\:dS.
$$
Since $U=$const. on $\p \Omega$,
\begin{align}\label{iso}
{\cal{E}}(\Omega,\rho) \leq {\cal{E}}(B_{r_\Omega},M).
\end{align}
\medskip

(ii) $G(s)=-H(s)$, where $H:\mathbb{R}^+\to \mathbb{R}^+$ and $H'(s)=h(s)>0.$

\noindent The comparison problem will be in this case
$$
\Delta \phi = \gamma^{n}h(\phi) \tx{in} B_{r_\Omega}, \quad \p_\nu \phi=\mu M \tx{on} \p B_{r_\Omega}.
$$
If $\mu M>0$ then  $\phi$ is increasing. Moreover we assume that $h$ is such that $\phi$ is positive. Under these  assumptions Lemma 3 yields for the transplanted function $U$ of $\phi$ 
$$
\int_\Omega H(U)\:dx \leq \gamma^n \int_{B_{r_\Omega}}H(\phi)\:dx.
$$
Consequently
$$
{\cal{E}}(\Omega)\leq \int_{B_{r_\Omega}} |\nabla \phi|^2\:dx +2\gamma^n\int_{B_{r_\Omega}} H(\phi)\:dx -2\phi(r_\Omega) \mu M.
$$
The right-hand side is the energy of the comparison problem. An example is $h(s)=c^2s$.
\bigskip

{\bf Acknowledgement} This paper was written during a visit at the Newton Institute in Cambridge. Both authors would like to thank this Institute for the excellent working atmosphere. They are also indebted to A. Henrot and D. Bucur who drove their attention to the problem concerning the principal eigenvalue. They thank in particular D. Bucur for having pointed out an error in the first draft of this note.


\end{document}